\documentclass[11pt]{amsart}

\usepackage[TS1,T1]{fontenc}
\usepackage[latin1]{inputenc}
\usepackage{amsfonts,amssymb,amsmath,amsthm,enumerate,latexsym}

\def\R{{\mathbb {R}}}
\def\N{{\mathbb {N}}}

\def\A{{\mathcal {A}}}
\def\eps{{\epsilon}}

\newtheorem{teo}{Theorem}[section]
\newtheorem{lema}[teo]{Lemma}
\newtheorem{prop}[teo]{Proposition}

\theoremstyle{remark}
\newtheorem{remark}[teo]{Remark}

\theoremstyle{definition}

\numberwithin{equation}{section}

\begin{document}

\title[Existence of solution]{Existence of solution to a critical equation with variable exponent}
\author[J. Fern\'andez Bonder, N. Saintier and A. Silva]{Juli\'an Fern\'andez Bonder, Nicolas Saintier and Anal\'ia Silva}

\address[J. Fern\'andez Bonder and A. Silva]{IMAS - CONICET and Departamento de Matem\'atica, FCEyN - Universidad de Buenos Aires, Ciudad Universitaria, Pabell\'on I  (1428) Buenos Aires, Argentina.}

\address[N. Saintier]{Instituto de Ciencias, Universidad Nacional de General Sarmiento, Juan María Gutierrez 1150 Los Polvorines - Pcia de Bs. As. - Argentina and Departamento de Matem\'atica, FCEyN - Universidad de Buenos Aires, Ciudad Universitaria, Pabell\'on I  (1428) Buenos Aires, Argentina.}

\email[J. Fernandez Bonder]{jfbonder@dm.uba.ar}

\urladdr[J. Fernandez Bonder]{http://mate.dm.uba.ar/~jfbonder}

\email[A. Silva]{asilva@dm.uba.ar}

\email[N. Saintier]{nsaintie@dm.uba.ar, nsaintie@ungs.edu.ar}

\urladdr[N. Saintier]{http://mate.dm.uba.ar/~nsaintie}


\subjclass[2010]{35J92,35B33}

\keywords{Sobolev embedding, variable exponents, critical exponents, concentration compactness}

\begin{abstract}
In this paper we study the existence problem for the $p(x)-$Laplacian operator with a nonlinear critical source. We find a local condition on the exponents ensuring the existence of a nontrivial solution that shows that the Pohozaev obstruction does not holds in general in the variable exponent setting. The proof relies on the Concentration--Compactness Principle for variable exponents and the Mountain Pass Theorem.
\end{abstract}

\maketitle

\section{Introduction}
In this paper we address the existence problem for the $p(x)-$Laplace operator with a source that has critical growth in the sense of the Sobolev embeddings. To be precise, we consider the equation
\begin{equation}\label{MainEq}
\begin{cases}
 -\Delta_{p(\cdot)}u + h|u|^{p(\cdot)-2}u  = |u|^{q(\cdot)-2}u &\quad \text{in } U, \\
 u = 0 &\quad \text{on } \partial U,
\end{cases}
\end{equation}
where $U\subset \R^n$ is a {\bf smooth} bounded domain, $p,q\colon U\to [1,\infty)$ are Log-H\"older continuous functions such that $1<\inf_U p\le \sup_U p<n$ and $1\le q(x)\le p^*(x):=np(x)/(n-p(x))$, $x\in U$.

The $p(x)$-Laplacian operator $\Delta_{p(\cdot)}$ is defined, as usual, as
$$
\Delta_{p(\cdot)} u := div(|\nabla u|^{p(\cdot)-2}\nabla u).
$$
This operator appears in the study of the so-called electrorheological fluids. We refer to the monograph by M. Ru{\v{z}}i{\v{c}}ka, \cite{Ru}, and its references, for a detailed account. In particular, after some simplifications, the modelling of these fluids lead to solve
\begin{equation}\label{ecuacion}
\begin{cases}
-\Delta_{p(\cdot)} u = f(\cdot,u,\nabla u) & \mbox{in }U,\\
u = 0 & \mbox{on }\partial U,
\end{cases}
\end{equation}
for some nonlinear source $f$. In most cases, the source term is taken to be only dependent on $u$ and so, in order for the usual variational techniques to work, one needs a control on the growth of $f$ given by the Sobolev embedding. 

When the growth of $f$ is subcritical in the sense of the Sobolev embedding the existence of solution follows easily by applying standard procedures of the calculus of variations (see e.g. \cite{Cabada-Pouso, Dinu, FZ, Mihailescu, Mihailescu-Radulescu} and many others).  On the other hand, when the source term has critical growth, there are only a few results on the existence of solutions for \eqref{ecuacion} that we will review thoroughly later. Let us just notice for the moment that these results only provide global existence conditions. This strongly contrasts with the constant exponent case that has been widely studied since Aubin and Brezis--Nirenberg's seminal works \cite{Aubin, BN}, and for which it is generally possible to provide existence conditions that are local in the sense that they involve the behaviour of the coefficients of the equations (and possibly some relevant geometric quantities) only in a neighborhood of a point.  Our main purpose in this paper is to provide local existence conditions for the equation \eqref{MainEq}.

\bigskip

In order to study \eqref{MainEq} by means of variational methods, we introduce the functional $J:W^{1,p(\cdot)}_0(U)\to\mathbb{R}$ defined by  
\begin{equation}\label{DefJ} 
J(u):= \int_U \frac{1}{p(x)}\Big(|\nabla u|^{p(x)} + h(x)|u|^{p(x)}\Big)\, dx - \int_U \frac{1}{q(x)} |u|^{q(x)}\, dx.
\end{equation} 
This functional is naturally associated to (\ref{MainEq}) in the sense that a weak solution of (\ref{MainEq}) is a critical point of $J$. We refer to Section 2 for the definition and some elementary properties of variable exponent spaces.

We need to assume that the smooth function $h$ is such that the functional
\begin{equation}\label{DefI}
 I(u):=\int_U \frac{1}{p(x)}\left(|\nabla u|^{p(x)}+h(x)|u|^{p(x)}\right)\,dx 
\end{equation} 
is coercive in the sense that the norm
$$
\|u\|:= \inf\,\left\{\lambda>0\, \int_U \left|\frac{\nabla u + h(x)u(x)}{\lambda}\right|^{p(x)}\,dx \le 1\right\} $$
is equivalent to the usual norm of $W_0^{1,p(\cdot)}(U)$.

\bigskip

When $\inf_U q > \sup_U p$, it is easy to show that $J$ satisfies the geometric assumptions of the Mountain--Pass Theorem (cf. Section 4). Hence if we assume moreover that the exponent $q$ is {\em subcritical} in the sense that
\begin{equation}\label{subcritical}
\inf_U (p^* - q) > 0,
\end{equation}
which implies that the immersion $W^{1,p(\cdot)}_0(U)\hookrightarrow L^{q(\cdot)}(U)$ is compact, then $J$ satisfies the Palais--Smale condition, and the existence of a nontrivial solution to \eqref{MainEq} follows easily.

When \eqref{subcritical} is violated, the immersion $W^{1,p(\cdot)}_0(U)\hookrightarrow L^{q(\cdot)}(U)$ does not need to be compact and so the Palais--Smale condition may fail. The existence of a non-trivial solution to \eqref{MainEq} is then a non-trivial problem.

We denote by 
\begin{equation}\label{CritSet} 
\A:=\{x\in U\colon q(x)=p^*(x)\}
\end{equation} 
 the {\it critical set}. We will assume in this work that this critical set is nonempty.

In \cite{MOSS} the authors prove that if $\A$ is {\em small} and there exists a control on the rate of how $q$ reaches the critical value $p^*$, then the immersion  $W^{1,p(\cdot)}_0(U)\hookrightarrow L^{q(\cdot)}(U)$ remains compact, and so the usual techniques can be applied. When the immersion fails to be compact they prove that if the subcriticality set $U\setminus \A$ contains a sufficiently large ball, then \eqref{MainEq} with $h=0$ has a nonnegative solution.

In \cite{FBS1}, problem \eqref{MainEq} is studied with $h=0$ and with a subcritical perturbation. In this work the authors generalize the Concentration--Compactness Principle (CCP) of P.L. Lions to the variable exponent case and prove that if the subcritical perturbation is large enough on the critical set, the Palais--Smale condition is verified and so the existence of a nontrivial solution follows. See also  \cite{Fu} where similar results were obtained independently.

In \cite{Silva}, using the CCP of \cite{FBS1, Fu}, a multiplicity problem for \eqref{MainEq} with $h=0$ and a nonsymmetric subcritical perturbation is analyzed.

More recently, the authors in \cite{FBSS1} studied the best Sobolev constant $S(p(.),q(.),U)$ corresponding to the embedding $W^{1,p(\cdot)}_0(U)\hookrightarrow L^{q(\cdot)}(U)$, namely 
\begin{equation}\label{BestSob} 
 S(p(.),q(.),U) = \inf_{u\in W^{1,p(.)}_0(U)} \dfrac{\|\nabla u\|_{L^{p(.)}(U)}}{\|u\|_{L^{q(.)}(U)}}. 
\end{equation} 
Using a refinement of the CCP proved in \cite{FBS1}, they gave sufficient conditions for the existence of an extremal for $S(p(.),q(.),U)$, and so the existence of a solution to \eqref{MainEq} with $h=0$ follows.

The study of \eqref{MainEq} posed in the whole $\R^n$ is analyzed in \cite{Alves, Fu2}. In those works the authors studied the problem in the case where $p, q$ and $h$ are radial functions and give somewhat restrictive conditions to ensure the existence of a nontrivial radial solution.

From now on we will assume that 
\begin{equation}\label{Hyp}  
 \sup_U p < \inf_U q. 
\end{equation} 

Our first result provides a condition for the functional $J$ defined by (\ref{DefJ}) to satisfy the Palais-Smale condition.

\begin{teo}\label{teoPScond}
The functional $J$ satisfies the Palais-Smale condition at level $c\in (0,\frac{1}{n}S^n)$ where
\begin{equation}\label{LocBestSob}
 S := \inf_{x\in\mathcal{A}} \lim_{\epsilon\to 0} S(p(\cdot),q(\cdot),B_\epsilon(x)),
\end{equation} 
and $S(p(\cdot),q(\cdot),B_\eps(x))$ stands for the best Sobolev constant for the domain $B_\eps(x)$ defined in a similar way as in (\ref{BestSob}). 
\end{teo}

\noindent The proof of Theorem \ref{teoPScond}  relies on a precise computation of the constants in the CCP recently proved by \cite{FBSS1}.

\medskip 

As a corollary, we can apply the Mountain--Pass Theorem to obtain the following necessary existence condition:

\begin{teo}\label{teoMP}
If there exists $v\in W_0^{1,p(\cdot)}(U)$ such that
\begin{equation}\label{CCPCond}
 \sup_{t>0} J(tv) < \frac{1}{n}S^n,
\end{equation}
then \eqref{MainEq} has a non-trivial nonnegative solution.
\end{teo}

Eventually the following result provide a sufficient local condition for (\ref{CCPCond}) to hold:

\begin{teo}\label{teocuentas}
Assume that the infimum in the definition \eqref{LocBestSob} of $S$ is attained at a point $x_0\in\mathcal{A}$ such that $x_0$ is a local minimum of $p$ and  a local maximum of $q$. In particular
\begin{equation}\label{des}
-\Delta p(x_0)\le 0\le - \Delta q(x_0).
\end{equation}
Assume moreover that $p,q$ are $C^2$ in a neighborhood of $x_0$, and that $h(x_0)<0$ if $1<p(x_0)<2$ ($n\ge 4$), or if $2\le p(x_0)<\sqrt{n}$ ($n\ge 5$), that at least one of the two inequalities in \eqref{des} is strict, but $h(x_0)$ is arbitrary. Under these assumptions \eqref{CCPCond} holds. In particular \eqref{MainEq} has a non-trivial nonnegative solution.
\end{teo}

 In the constant exponent case, the well known Pohozaev obstruction \cite{Pohozaev} affirms that if $h\ge 0$ and $U$ is starshaped then there are no (positive) solutions to \eqref{MainEq}. Our result shows that for variable $p$ and $q$ and $p(x)\ge 2$ this does not need to be the case, showing a stricking difference between the constant exponent case and the variable exponent one.

\section{Preliminaries on variable exponent spaces.}

In this section we review some preliminary results regarding Lebesgue and Sobolev spaces with variable exponent. All of these results and a comprehensive study of these spaces can be found in \cite{libro}. 

Consider a function $p:U\to [1,+\infty]$ Log-H\"older continuous in the sense that 
$$
|p(x)-p(y)| \le \frac{C}{|\log|x-y||},\quad \text{for } x,y\in U,\ x\neq y 
$$
for some constant $C>0$. This regularity assumptions is not needed to define the Lebesgue and Sobolev spaces with variable $p$ but turns out to be very useful for these Sobolev spaces to enjoy all the usual properties like Sobolev embeddings, Poincar\'e inequality and so on. We will therefore assume it from now for simplicity.

The variable exponent Lebesgue space $L^{p(x)}(U)$ is defined by
$$
L^{p(x)}(U) = \Big\{u\in L^1_{\text{loc}}(U) \colon \int_U|u(x)|^{p(x)}\,dx<\infty\Big\}.
$$
This space is endowed with the norm
$$
\|u\|_{L^{p(x)}(U)}=\inf\Big\{\lambda>0:\int_U\Big|\frac{u(x)}{\lambda}\Big|^{p(x)}\,dx\leq 1\Big\}.
$$
The variable exponent Sobolev space $W^{1,p(x)}(U)$ is defined by
$$
W^{1,p(x)}(U) = \{u\in W^{1,1}_{\text{loc}}(U) \colon u\in L^{p(x)}(U) \mbox{ and } |\nabla u |\in  L^{p(x)}(U)\}.
$$
The corresponding norm for this space is
$$
\|u\|_{W^{1,p(x)}(U)}=\|u\|_{L^{p(x)}(U)}+\| \nabla u \|_{L^{p(x)}(U)}.
$$
Define $W^{1,p(x)}_0(U)$ as the closure of $C_c^\infty(U)$ with respect to the $W^{1,p(x)}(U)$ norm. The spaces $L^{p(x)}(U)$, $W^{1,p(x)}(U)$ and $W^{1,p(x)}_0(U)$ are separable and reflexive Banach spaces when 
$1<p^- \le p^+ <\infty$, where $p^-:=ess-\inf_U p$ and $p^+:=ess-\sup_U p$.

As usual, we denote the conjugate exponent of $p(x)$ by $p'(x) = p(x)/(p(x)-1)$  and the Sobolev exponent by
$$
p^*(x)=\begin{cases}
\frac{Np(x)}{N-p(x)} & \mbox{ if } p(x)<N,\\
\infty & \mbox{ if } p(x)\geq N.
\end{cases}
$$

The following result is proved in \cite{Fan} (see also \cite{libro}, pp. 79, Lemma 3.2.20 (3.2.23)).

\begin{prop}[H\"older-type inequality]\label{Holder}
Let $f\in L^{p(x)}(U)$ and $g\in L^{q(x)}(U)$. Then the following inequality holds
$$
\| fg \|_{L^{s(x)}(U)}\le  \Big( \Big(\frac{s}{p}\Big)^+ + \Big(\frac{s}{q}\Big)^+\Big) \|f\|_{L^{p(x)}(U)}\|g\|_{L^{q(x)}(U)},
$$
where
$$
\frac{1}{s(x)} = \frac{1}{p(x)} + \frac{1}{q(x)}.
$$
\end{prop}

The Sobolev embedding Theorem is also proved in \cite{Fan}, Theorem 2.3.

\begin{prop}[Sobolev embedding]\label{embedding}
Let $q:U\in [1,+\infty)$ be a measurable function such that $1\leq q(x)\le p^*(x)<\infty$ for all $x\in\overline{U}$.  Then there is a continuous embedding
$$ W^{1,p(x)}(U)\hookrightarrow L^{q(x)}(U). $$
Moreover, if $\inf_{U} (p^*-q)>0$ then, the embedding is compact.
\end{prop}

As in the constant exponent spaces, Poincar\'e inequality holds true (see \cite{libro}, pp. 249, Theorem 8.2.4)
\begin{prop}[Poincar\'e inequality]\label{Poincare}
There is a constant $C>0$, $C=C(U)$, such that
$$ \|u\|_{L^{p(x)}(U)}\leq C\|\nabla u\|_{W^{1,p(x)}(U)}, $$
for all $u\in W^{1,p(x)}_0(U)$.
\end{prop}

\noindent It follows in particular from the Poincar\'e inequality that $\| \nabla u \|_{L^{p(x)}(U)}$ and $\|u\|_{W^{1,p(x)}(U)}$ are equivalent norms on $W_0^{1,p(x)}(U)$.

Throughout this paper the following notation will be used: Given $q\colon U\to\R$ bounded, we denote
$$
q^+ := \sup_U q(x), \qquad q^- := \inf_U q(x).
$$

The following proposition is also proved in \cite{Fan} and it will be most usefull (see also \cite{libro}, Chapter 2, Section 1).
\begin{prop}\label{norma.y.rho}
Set $\rho(u):=\int_U|u(x)|^{p(x)}\,dx$. For $u,\in L^{p(x)}(U)$ and $\{u_k\}_{k\in\N}\subset L^{p(x)}(U)$, we have
\begin{align}
& u\neq 0 \Rightarrow \Big(\|u\|_{L^{p(x)}(U)} = \lambda \Leftrightarrow \rho(\frac{u}{\lambda})=1\Big).\\
& \|u\|_{L^{p(x)}(U)}<1 (=1; >1) \Leftrightarrow \rho(u)<1(=1;>1).\\
& \|u\|_{L^{p(x)}(U)}>1 \Rightarrow \|u\|^{p^-}_{L^{p(x)}(U)} \leq \rho(u) \leq \|u\|^{p^+}_{L^{p(x)}(U)}.\\
& \|u\|_{L^{p(x)}(U)}<1 \Rightarrow \|u\|^{p^+}_{L^{p(x)}(U)} \leq \rho(u) \leq \|u\|^{p^-}_{L^{p(x)}(U)}.\\
& \lim_{k\to\infty}\|u_k\|_{L^{p(x)}(U)} = 0 \Leftrightarrow \lim_{k\to\infty}\rho(u_k)=0.\\
& \lim_{k\to\infty}\|u_k\|_{L^{p(x)}(U)} = \infty \Leftrightarrow \lim_{k\to\infty}\rho(u_k) = \infty.
\end{align}
\end{prop}

The following Lemma is the extension to variable exponents of the well-known Brezis-Lieb Lemma (see \cite{Brezis-Lieb}). The proof is analogous to that of \cite{Brezis-Lieb}. See Lemma 3.4 in \cite{FBS1}. 

\begin{lema}\label{Brezis-Lieb}
Let $f_n\to f$ a.e and $f_n\rightharpoonup f$ in $L^{p(x)}(U)$ then
$$
\lim_{n\to\infty}\left(\int_U|f_n|^{p(x)}dx-\int_U|f-f_n|^{p(x)}dx\right)=\int_U|f|^{p(x)}dx.
$$
\end{lema}
For much more on these spaces, we refer to \cite{libro}.

\section{Proof of theorem \ref{teoPScond}}

In this section we verify that the functional $J$ defined by (\ref{DefJ}) satisfies the Palais--Smale condition (PS for short) for energy levels below the critical one $\tfrac{1}{n} S^n$. The scheme of the proof is classical (see e.g. \cite{Saintier}) but relies on a version of Lions' concentration--compactness principle adapted to the variable exponent setting in \cite{FBS1} and then refined in \cite{FBSS1}. 

\medskip

Let $\{u_k\}_{k\in\N}\subset W^{1,p(\cdot)}(U)$ be a PS--sequence for $J$. Recall that this means that the sequence $\{J(u_k)\}_{k\in\N}$ is bounded, and that $DJ(u_k)\to 0$ strongly in the dual space $W^{1,p(.)}(U)'$.

Recalling that the functional $I$ defined by (\ref{DefI}) is assumed to be coercive, it then follows that $\{u_k\}_{k\in\N}$ is bounded in $W^{1,p(\cdot)}(U)$. In fact, for $k$ large, we have that
\begin{align*}
c+1 &\ge J(u_k) - \frac{1}{q^-}\langle DJ(u_k), u_k\rangle \\
&\ge \big(\frac{1}{p^+}-\frac{1}{q^-}\big) \int_U |\nabla u_k|^{p(x)} + h(x) |u_k|^{p(x)}\, dx -
\int_U \big(\frac{1}{q(x)}-\frac{1}{q^-}\big) |u_k|^{q(x)}\, dx\\
&\ge \big(\frac{1}{p^+}-\frac{1}{q^-}\big) \int_U |\nabla u_k|^{p(x)} + h(x) |u_k|^{p(x)}\, dx.
\end{align*}
from where the claim follows recalling assumption (\ref{Hyp}).

We may thus assume that $u_k\rightharpoonup u$ weakly in $W_0^{1,p(\cdot)}(U)$. We claim that $u$ turns out to be a weak solution to \eqref{MainEq}. In fact, since $u_k\rightharpoonup u$ weakly in $W^{1,p(\cdot)}_0(U)$ it follows that
\begin{align*}
&|\nabla u_k|^{p(\cdot)-2}\nabla u_k \rightharpoonup |\nabla u|^{p(\cdot)-2}\nabla u \quad\mbox{weakly in }L^{p'(\cdot)}(U),\\
&|u_k|^{p(\cdot)-2} u_k \rightharpoonup |u|^{p(\cdot)-2} u \quad\mbox{weakly in }L^{p'(\cdot)}(U),\\
&|u_k|^{q(\cdot)-2} u_k \rightharpoonup |u|^{q(\cdot)-2} u \quad\mbox{weakly in }L^{q'(\cdot)}(U).
\end{align*}
So
\begin{equation}\label{conv.p}
\begin{split}
0&=\lim_{k\to\infty}\langle DJ(u_k), \phi\rangle\\
&= \lim_{k\to\infty} \int_U |\nabla u_k|^{p(x)-2}\nabla u_k \nabla \phi + h |u_k|^{p(x)-2} u_k\phi\, dx - \int_U |u_k|^{q(x)-2} u_k \phi\\
&= \int_U |\nabla u|^{p(x)-2}\nabla u \nabla \phi + h |u|^{p(x)-2} u\phi\, dx - \int_U |u|^{q(x)-2} u \phi,
\end{split}
\end{equation}
for every $\phi\in C^\infty_0(U)$. This proves that $u$ is a weak solution of \eqref{MainEq}. 

By the CCP for variable exponents (see \cite{FBS1} and the refinement proved in \cite{FBSS1}) it holds that
\begin{align*}
& |u_k|^{q(\cdot)} \rightharpoonup \nu = |u|^{q(\cdot)} + \sum_{i\in I} \nu_i \delta_{x_i} \quad \mbox{weakly in the sense of measures,}\\
& |\nabla u_k|^{p(\cdot)} \rightharpoonup \mu \ge |\nabla u|^{p(\cdot)} + \sum_{i\in I} \mu_i \delta_{x_i} \quad \mbox{weakly in the sense of measures,}\\
& S\nu_i^{1/p^*(x_i)}\le \mu_i^{1/p(x_i)},
\end{align*}
where $I$ is a finite set, $\{\nu_i\}_{i\in I}$ and $\{\mu_i\}_{i\in I}$ are positive numbers and the points $\{x_i\}_{i\in I}$ belong to the critical set $\A$ defined in (\ref{CritSet}).

It is not difficult to check (arguing as in \eqref{conv.p}) that $v_k:=u_k-u$ is PS--sequence for $\tilde J(v):=J(v)-\int_U \frac{1}{p(x)} h|v|^{p(x)}$. Now, by Lemma \ref{Brezis-Lieb} we get
\begin{align*}
J(u_k)-J(u) & =   \int_U \frac{1}{p(x)}\Big[|\nabla v_k|^{p(x)} + h |v_k|^{p(x)}\Big]\, dx - \int_U \frac{1}{q(x)} |v_k|^{q(x)}\, dx + o(1)\\
& =  \tilde J(v_k) + \int_U  \frac{1}{p(x)}  h |v_k|^{p(x)}\, dx + o(1)\\
& =  \tilde J(v_k) + o(1).
\end{align*}

Since $u$ is a weak solution of \eqref{MainEq},  and since $p^+<q^-$,
\begin{equation*}
\begin{split}
 J(u) & \ge \frac{1}{p^+} \int_U \left(|\nabla u|^{p(x)}+h(x)|u|^{p(x)}\right)\,dx
           - \frac{1}{q^-} \int_U |u|^{q(x)}\,dx  \\
      & = \left( \frac{1}{p^+} - \frac{1}{q^-} \right) \int_U |u|^{q(x)}\,dx \\
      & \ge 0.
\end{split}
\end{equation*}
Therefore,
$$
J(u_k)\ge \tilde J(v_k) + o(1).
$$

Let $\phi\in C^\infty_c(U)$. As $D\tilde J(v_k)\to 0$, we have
\begin{align*}
o(1) &= \langle D\tilde J(v_k), v_k\phi\rangle \\
&= \int_U |\nabla v_k|^{p(x)}\phi\, dx - \int_U |v_k|^{q(x)}\phi\, dx + \int_U |\nabla v_k|^{p(x)-2}\nabla v_k \nabla\phi v_k\, dx\\
&= A - B + C.
\end{align*}
Since $v_k\rightharpoonup 0$ weakly in $W^{1,p(\cdot)}_0(U)$ it is easy to see that $C\to 0$ as $k\to\infty$. By means of Lemma \ref{Brezis-Lieb} it follows that
$$
A\to \int_U\phi\, d\tilde\mu \qquad \mbox{and}\qquad B\to \int_U \phi\, d\tilde\nu,
$$
where $\tilde\mu = \mu - |\nabla u|^{p(x)}$ and $\tilde\nu =\nu - |u|^{q(x)}$. So we conclude that $\tilde\mu=\tilde\nu$. In particular $\nu_i\ge \mu_i$ ($i\in I$) from where we obtain that $\nu_i\ge S^n$. Hence
\begin{equation*}
\begin{split}
 c &=  \lim_{k\to\infty} J(u_k) \ge \lim_{k\to\infty} \tilde J(v_k)  =\int \frac{1}{p(x)}\,d\tilde\mu - \int \frac{1}{q(x)}\,d\tilde\nu\\
 & = \int \Big(\frac{1}{p(x)}-\frac{1}{q(x)}\Big) \,d\tilde\nu = \sum_{i\in I} \left(\frac{1}{p(x_i)}-\frac{1}{p^*(x_i)}\right) \nu_i \\
 & \ge \#(I)\frac{1}{n} S^n.
\end{split}
\end{equation*}
We deduce that if $c<\frac{1}{n} S^n$ then $I$ must be empty implying that $u_k\to u$ strongly in $W^{1,p(\cdot)}(U)$.

\section{Proof of theorem \ref{teoMP}}

The proof of Theorem \ref{teoMP} is an immediate consequence of the Mountain--Pass Theorem, Theorem \ref{teoPScond} and assumption \eqref{CCPCond}.

In fact, it suffices to verify that $J$ has the Mountain--Pass geometry and that $J(tu)<0$ for some $t>0$. Concerning the latter condition notice that for $t>1$,
\begin{equation*}
\begin{split}
 J(tu) & = \int_U \frac{t^{p(x)}}{p(x)}\left(|\nabla u|^{p(x)}+h(x)|u|^{p(x)}\right)\,dx
           - \int_U \frac{t^{q(x)}}{q(x)}|\nabla u|^{q(x)}\,dx \\
       & \le t^{p^+}I(u) - t^{q^-}\int_U \frac{1}{q(x)}|\nabla u|^{q(x)}\,dx,
\end{split}
\end{equation*}
which tends to $-\infty$ as $t\to +\infty$ since $q^->p^+$.

It remains to see that $J$ has the Mountain--Pass geometry. But $J(0)=0$ and, if $\|v\|_{W_0^{1,p(\cdot)}(U)}=r$ small enough, then
$$
\int_U |\nabla v|^{p(x)} + h |v|^{p(x)}\, dx \ge {\bf C} \| v\|_{W_0^{1,p(\cdot)}(U)}^{p^+}
$$
and
$$
\|v\|_{L^{q(\cdot)}(U)} \le C\|v\|_{W_0^{1,p(\cdot)}(U)} = Cr<1,
$$
so
$$
\int_U |v|^{q(x)}\, dx \le C\|v\|_{W_0^{1,p(\cdot)}(U)}^{q^-}.
$$
Therefore
$$
J(v) \ge \frac{{\bf C}}{p^+} r^{p^+} - \frac{C}{q^-} r^{q^-} >0,
$$
since $p^+<q^-$. This completes the proof.

\section{Proof of theorem \ref{teocuentas}}

Let $x_0\in \mathcal{A}$ be such that
$$
 S := \inf_{x\in\mathcal{A}} \lim_{\epsilon\to 0} S(p(\cdot),q(\cdot),B_\epsilon(x)) = \lim_{\epsilon\to 0} S(p(\cdot),q(\cdot),B_\epsilon(x_0)).
$$

For ease of notation we assume that $x_0=0$, write $p=p(0)$ and observe that $q=q(0) = p^*$.
From Theorem 6.1 in \cite{FBSS1}, we have that if $0$ is a local maximum of $q$ and a local minimum of $p$, then 
$$S=\lim_{\eps\to 0}S(p(\cdot),q(\cdot),B_\eps(0))=K(n,p)^{-1},$$
where $K(n,p)$ is the best constant in the Sobolev inequality in $\R^n$, i.e.
$$
K(n,p)^{-1} = \inf_{v\in W^{1,p}(\R^n)} \frac{\|\nabla v\|_{L^p(\R^n)}}{\| v\|_{L^{p^*}(\R^n)}}.
$$
Let $U$ be an extremal for the constant $K(n,p)$. That is, $U$ verifies
$$
K(n,p)^{-1}  = \frac{\|\nabla U\|_{L^p(\R^n)}}{\| U\|_{L^{p^*}(\R^n)}}.
$$

It is well known, see \cite{Aubin, Talenti}, that $U$ can be given by the formula
$$ U(x)=\left(1+|x|^\frac{p}{p-1}\right)^{-\frac{n-p}{p}}. $$
Moreover, any extremal for $K(n,p)$ is obtained by a translation and a dilation of $U$ in the form
$$
U_{\eps, x_0}(x) = \eps^{-\frac{n-p}{p}} U((x-x_0)/\eps).
$$

Given $\delta>0$ small we take a cut-off function $\eta\in C^\infty_c(B_{2\delta},[0,1])$ such that $\eta\equiv 1$ in $B_\delta$. We then consider the test-function
$$ u_\eps(x) =  U_{\eps,0}(x)\eta(x). $$
For this test function we have:
\begin{prop}\label{PropCuentas}
Assume that $0$ is a critical point of $p$ and $q$. We have
\begin{itemize}
\item If $p\le \frac{n}{2}$,
 \begin{equation}\label{EstimLq}
  \int_{\R^n} f(x)u_\eps^{q(x)}\,dx =
   A_0 + A_1 \eps^2\ln\eps + o(\eps^2\ln\eps)
 \end{equation}
 with
 \begin{equation*}
  A_0 = f(0)\int_{\R^n} U^{p^*}\,dx, \quad
  A_1 = - \frac{n-p}{p} \frac{f(0)}{2} \int_{\R^n} U^{p^*} (D^2q(0)x,x) \,dx
 \end{equation*}

\item  If $p<\min\{\sqrt{n},\frac{n+2}{3}\}$,
 \begin{equation}\label{EstimGrad}
  \int_{\R^n} f(x)|\nabla u_\eps|^{p(x)}\,dx =
  B_0 +  B_1 \eps^2\ln\eps + o(\eps^2\ln\eps)
 \end{equation}
 with
 \begin{equation*}
  B_0 = f(0)\int_{\R^n} |\nabla U|^p\,dx, \qquad
  B_1 = - \frac{n}{p} \frac{f(0)}{2} \int_{\R^n} |\nabla U|^p (D^2p(0)x,x) \,dx
 \end{equation*}

\item If $p<\sqrt{n}$,
 \begin{equation}\label{EstimLp}
  \int_{\R^n} f(x)|u_\eps|^{p(x)}\,dx =  C_0 \eps^p + o(\eps^p)
  \quad \text{with } \quad C_0  = f(0)\int_{\R^n} U^p\,dx.
 \end{equation}
 \end{itemize}
\end{prop}

\begin{remark}
Observe that if $g(x)$ is a radial function then
$$ \int_{\R^n} g(x) (Ax,x)\, dx = tr(A) \int_{\R^n} g(x) x_1^2\, dx
                                = \dfrac{tr(A)}{n} \int_{\R^n} g(x) |x|^2\, dx,
$$
for any $A\in \R^{n\times n}$ (with adequate decaying assumptions at infinity on $g$).
In fact this is a consequence of the fact that, for $i\neq j$,
$$ \int_{\R^n} g(x) x_i x_j\, dx = 0. $$
With this observation, we easily conclude that
$$
A_1 = - \frac{f(0)}{p^*} \Delta q(0) \int_{\R^n} U^{p^*} |x|^2 \,dx
$$
and
$$
B_1 = - \frac{f(0)}{2p} \Delta p(0) \int_{\R^n} |\nabla U|^p |x|^2 \,dx.
$$
\end{remark}

We postpone the proof of this proposition to Section 6.

\bigskip 

As $U$ is an extremal for $K(n,p)$ it follows that $U$ verifies
$$ 
-\Delta_p U = \frac{K(n,p)^{-p}}{\|U\|_{L^{p^*}(\R^n)}^{p^*-p}} U^{p^*-1} = CU^{p^*-1}. 
$$
Then $V=C^\frac{1}{p^*-p}U = \frac{K(n,p)^{-\frac{n-p}{p}}}{\|U\|_{p^*}}U$ solves $-\Delta_p V=V^{p^*-1}$ and satisfy 
$$
\|\nabla V\|_{L^p(\R^n)} = K(n,p)^{-n/p}.
$$

Consider the test function
$$ v_\eps(x) = \eps^{-\frac{n-p}{p}} V(x/\eps)\eta(x) = C^\frac{1}{p^*-p} u_\eps(x). $$
Using the previous proposition we immediately obtain

\begin{prop} Assume that $0$ is a critical point of $p$ and $q$. If $p<\min\{\sqrt{n},\frac{n+2}{3}\}$ then
 \begin{equation}\label{EstimLqbis}
 \begin{split}
  \int_{\R^n} f(x)v_\eps^{q(x)}\,dx & =
   f(0)K(n,p)^{-n}+ f(0) A \eps^2\ln\eps + o(\eps^2\ln\eps), \\
  \int_{\R^n} f(x)|\nabla v_\eps|^{p(x)}\,dx & =
  f(0)K(n,p)^{-n} +  f(0) B \eps^2\ln\eps + o(\eps^2\ln\eps), \\
  \int_{\R^n} f(x)|v_\eps|^{p(x)}\,dx & = f(0) C \eps^p + o(\eps^p),
 \end{split}
 \end{equation}
 with
 \begin{equation*}
 \begin{split}
   A & = -\frac{\Delta q(0)}{2p^*}K(n,p)^{-n} \|U\|_{p^*}^{-p^*} \int_{\R^n} |x|^2U^{p^*}\,dx,\\
   B & =  -\frac{\Delta p(0)}{2p}K(n,p)^{p-n}\|U\|_{p^*}^{-p} \int_{\R^n} |x|^2|\nabla U|^p \,dx,\\
   C & = K(n,p)^{p-n}\|U\|_{p^*}^{-p} \|U\|_p^p.
 \end{split}
 \end{equation*}
\end{prop}

Using $v_\eps$ as a test-function in (\ref{CCPCond}) we can see that there exists $t_0>1$ such that $J(tv_\eps)<0$ for $t>t_0$. Now if $p<2$, we can write
$$ 
f_\eps(t):=J(tv_\eps) = f_0(t) + \eps^p f_1(t) + o(\eps^p) 
$$
$C^1$-uniformly in $t\in [0,t_0]$, with
$$ 
f_0(t) = K(n,p)^{-n} \left(\frac{t^p}{p} - \frac{t^{p^*}}{p^*}\right), \quad \text{and} \quad
f_1(t) = \frac{1}{p}t^p h(0)C.
$$
Notice that $f_0$ reaches its maximum in $[0,t_0]$ at $t=1$. Moreover it is a nondegenerate maximum since $f_0''(1)=(p-p^*)K^{-n}\neq 0$. It follows that $f_\eps$ reaches a maximum at $t_\eps = 1 + a \eps^p + o(\eps^p)$ for $a=-\frac{f_1'(1)}{f_0''(1)}$. Hence
\begin{equation*}
\begin{split}
 \sup_{t>0}J(tv_\eps) & =J(t_\eps v_\eps) = \frac{1}{n}K(n,p)^{-n} + f_1(1)\eps^p + o(\eps^p)
\end{split}
\end{equation*}
Then if $h(0)<0$ we get $\sup_{t>0}J(tv_\eps)<\frac{1}{n}K(n,p)^{-n}$.

We now assume that $p\ge 2$. Then
\begin{equation*}
f_\eps(t)=J(tv_\eps) = f_0(t) + \tilde f_1(t) \eps^2\ln\eps  + o(\eps^2\ln\eps),
\end{equation*}
$C^1$-uniformly in $t\in [0,t_0]$, with
\begin{equation*}
 \tilde f_1(t) = \frac{t^{p^*}}{p^*} A -\frac{t^p}{p} B.
\end{equation*}
As before $f_\eps$ reaches its maximum at $t_\eps= 1+ a\eps^2\ln\eps + o(\eps^2\ln\eps)$ with $a=-\frac{\tilde f_1'(1)}{f_0''(1)}$. Hence
\begin{equation*}
\begin{split}
 \sup_{t>0}J(tv_\eps) & =J(t_\eps v_\eps) = f_0(1)+  \tilde f_1(1)\eps^2\ln\eps + o(\eps^2\ln\eps) \\
                      & = \frac{1}{n}K(n,p)^{-n} +  \tilde f_1(1)\eps^2\ln\eps + o(\eps^2\ln\eps).
\end{split}
\end{equation*}
We thus need $\tilde f_1(1)<0$ i.e.
\begin{equation}\label{cond}
 -\Delta p(0) < -\Delta q(0) (p/p^*)^2 D(n,p), \quad \text{where}\quad
         D(n,p):= \frac{\displaystyle \int_{\R^n} |\nabla U|^p\, dx  \int_{\R^n} |x|^2U^{p^*}\, dx}
                 {\displaystyle \int_{\R^n} U^{p^*}\, dx  \int_{\R^n} |x|^2|\nabla U|^p\, dx}.
\end{equation}
Since $0$ is a local maximum of $q$ and a local minimum of $p$ we already know that (\ref{des}) holds. Then if one of the two inequalities in (\ref{des}) is strict we see that (\ref{cond}) holds.

This ends the proof of Theorem \ref{teocuentas}.

\medskip

As a final remark, we notice that we can compute $D(n,p)$ exactly. To do this let
\begin{equation}\label{Ipq}
 I_p^q := \int_0^\infty t^{q-1}(1+t)^{-p}\,dt
  = B(q,p-q)
  = \frac{\Gamma (q)\Gamma (p-q)}{\Gamma (p)},
\end{equation}
where $B(x,y):=\int_0^\infty t^{x-1}(1+t)^{-x-y}\,dt $ is the Beta function. This formula can be found, for instance, in \cite{beta}.
Passing to spherical coordinates and then performing the change of variable $t=r^\frac{p}{p-1}$, $dr=\frac{p-1}{p}t^{-\frac{1}{p}}dt$, we obtain
\begin{equation*}
\begin{split}
 & \int_{\R^n} U^{p^*}\, dx = U_{n-1}\frac{p-1}{p} I_n^{n\frac{p-1}{p}}, \\
&  \int_{\R^n} |x|^2 U^{p^*}\, dx = \omega_{n-1}\frac{p-1}{p} I_n^{n\frac{p-1}{p} -\frac{2}{p} +2}, \\
 & \int_{\R^n} |\nabla U|^p\, dx = \omega_{n-1}\frac{p-1}{p} \left(\frac{n-p}{p-1}\right)^p I_n^{n\frac{p-1}{p} + 1}, \\
&  \int_{\R^n} |x|^2|\nabla U|^p\, dx = \omega_{n-1}\frac{p-1}{p} \left(\frac{n-p}{p-1}\right)^p I_n^{n\frac{p-1}{p} - \frac{2}{p} + 3}.
\end{split}
\end{equation*}
Then
\begin{equation*}
\begin{split}
 D(n,p) & = \frac{\displaystyle I_n^{\frac{n(p-1)}{p} + 1} I_n^{\frac{n(p-1)}{p} -\frac{2}{p} +2} }
               {\displaystyle I_n^{\frac{n(p-1)}{p} } I_n^{\frac{n(p-1)}{p} -\frac{2}{p} +3} }
         = \frac{n}{n-p}\frac{(n-p) - 2(p-1)}{n+2},
\end{split}
\end{equation*}
where we used that
$$ I_p^{q+1} = \frac{q}{p-q-1}I_p^q $$
which follows from (\ref{Ipq}) and the formula $\Gamma(z+1)=z\Gamma(z)$.

\section{Proof of Proposition \ref{PropCuentas}}

As $0$ is a local minimum of $p(\cdot)$ we can assume that $p^-_{2\delta}:=\min_{x\in B_{2\delta}}\,p(x) = p$.

\subsection{Proof of \eqref{EstimLq}}

We first write
$$ \int_{\R^n} f(x) u_\eps(x)^{q(x)}\,dx
   = \int_{B_{2\delta}\backslash B_{\eps^{1/p}}} f(x) u_\eps^{q(x)}\,dx
    + \int_{B_{\eps^{1/p}}} f(x) u_\eps(x)^{q(x)}\,dx
   = I_1(\eps) +   I_2(\eps).
$$
Since $u_\eps(x)\le 1$ if $|x|\ge \eps^{1/p}$, we have, letting $q^-_{2\delta}:=\min_{B_{2\delta}}\,q$ that
\begin{equation*}
\begin{split}
I_1(\eps) & \le \|f\|_{L^{\infty}(B_{2\delta})}   \int_{B_{2\delta}\backslash B_{\eps^{1/p}}} u_\eps(x)^{q^-_{2\delta}}\,dx  \\
& \le \|f\|_{L^{\infty}(B_{2\delta})}  \eps^{n-\frac{n-p}{p}q^-_{2\delta}}
        \int_{\R^n\backslash B_{\eps^{-(p-1)/p}}} U(x)^{q^-_{2\delta}}\,dx,
\end{split}
\end{equation*}
where the integral in the right hand side can be bounded by
\begin{equation*}
\begin{split}
  C \int_{\eps^{-(p-1)/p}}^{+\infty} (1+r^\frac{p}{p-1})^{-\frac{n-p}{p}q^-_{2\delta}} \, r^{n-1}\,dr
 & \le C  \int_{\eps^{-(p-1)/p}}^{+\infty} r^{-1+n-\frac{n-p}{p-1}q^-_{2\delta}}\,dr
   \le C  \eps^{-n\frac{p-1}{p}+\frac{n-p}{p}q^-_{2\delta}}.
\end{split}
\end{equation*}
Hence $I_1(\eps) \le C \eps^{n/p}$ so that
\begin{equation*}
\begin{split}
 \int_{\R^n} f(x) u_\eps(x)^{q(x)}\,dx
 & = \int_{B_{\eps^{1/p}}} f(x) u_\eps(x)^{q(x)}\,dx + O(\eps^{n/p}) \\
 & = \int_{B_{\eps^{-(p-1)/p}}} f(\eps x) \eps^{n-q(\eps x)\frac{n-p}{p}}U(x)^{q(\eps x)}\,dx + O(\eps^{n/p}).
\end{split}
\end{equation*}
As $\nabla q (0) = 0$ we get
$$ q(\eps x) = q(0) + \frac{1}{2}\eps^2 (D^2q(0)x,x) + o(\eps^2 |x|^2), $$
with $q(0)=p(0)^*=p^*$, so
\begin{equation*}
\begin{split}
 \int_{\R^n} f(x) u_\eps(x)^{q(x)}\,dx  = &  A_0(\eps) + A_1(\eps)\eps^2\ln\,\eps
   +    \int_{B_{\eps^{-(p-1)/p}}} o(\eps^2\ln\,\eps)|x|^2 U(x)^{p^*}\,dx\\
   & + \eps \int_{B_{\eps^{-(p-1)/p}}} U(x)^{p^*} \nabla f(0)\cdot x \, dx + O(\eps^{n/p}) \\
 = & A_0(\eps)  + A_1(\eps)\eps^2\ln\,\eps   +   o(\eps^2\ln\,\eps) + O(\eps^{n/p}),
\end{split}
\end{equation*}
where $A_0(\eps)$ and $A_1(\eps)$ are the same as $A_0$ and $A_1$ except that we integrate over $B_{\eps^{-(p-1)/p}}$ instead of $\R^n$ and we have used the fact that
$$
\int_{B_{\eps^{-(p-1)/p}}} U(x)^{p^*} \nabla f(0)\cdot x \, dx=0,
$$
since $U$ is radially symmetric. We have
\begin{equation*}
\begin{split}
 |A_0(\eps) - A_0|
 & \le C \int_{\R^n\backslash B_{\eps^{-(p-1)/p}}} U(x)^{p^*}\,dx\\
 & \le C \int_{\eps^{-(p-1)/p}}^{+\infty} (1+r^\frac{p}{p-1})^{-n}r^{n-1}\,dr \\
 & \le C \int_{\eps^{-(p-1)/p}}^{+\infty} r^{\frac{-np}{p-1} + n-1}\,dr \\
 & \le C  \eps^\frac{n}{p}.
\end{split}
\end{equation*}
If $p<(n+2)/2$, we can estimate
\begin{equation*}
\begin{split}
  |A_1(\eps)-A_1|
 & \le C \int_{\R^n\backslash B_{\eps^{-(p-1)/p}}} |x|^2U(x)^{p^*}\,dx\\
 &  \le C \int_{\eps^{-(p-1)/p}}^{+\infty} (1+r^\frac{p}{p-1})^{-n}r^{n+1}\,dr \\
 &  \le C \eps^\frac{n+2-2p}{p}.
\end{split}
\end{equation*}
We thus have
\begin{equation*}
\begin{split}
\int_{\R^n} f(x) u_\eps(x)^{q(x)}\,dx -  A_0  - A_1 \eps^2\ln\eps  
=  O(\eps^{n/p}) + o(\eps^2\ln\,\eps),
\end{split}
\end{equation*}
which reduces to (\ref{EstimLq}) if we assume that $p\le n/2$.

\subsection{Proof of \eqref{EstimLp}}

As before,
$$ \int_{\R^n} f(x) u_\eps^{p(x)}\,dx = \int_{B_{\eps^{1/p}}} f(x) u_\eps^{p(x)}\,dx +
                                  \int_{B_{2\delta}\backslash B_{\eps^{1/p}}} f(x) u_\eps^{p(x)}\,dx
$$
where, noticing that $p=p^-_{2\delta}$, the 2nd integral in the right hand side can be bounded by
\begin{equation*}
\begin{split}
 \int_{B_{2\delta}\backslash B_{\eps^{1/p}}} u_\eps^p\,dx
 & \le C \eps^p \int_{\eps^{1/p-1}}^\infty (1+r^\frac{p}{p-1})^{p-n}r^{n-1}\,dr\\
 & \le C \eps^p \eps^\frac{n-p^2}{p}\\
 & = C \eps^\frac{n}{p},
\end{split}
\end{equation*}
if $p^2<n$.
Then
\begin{equation*}
\begin{split}
  \int_{\R^n} f(x)u_\eps^{p(x)}\,dx
 & = \int_{B_{\eps^{1/p}}} f(x)u_\eps^{p(x)}\,dx + O(\eps^\frac{n}{p})\\
 & = \int_{B_{\eps^{1/p-1}}} f(\eps x) \eps^{n-\frac{n-p}{p}p(\eps x)} U(x)^{p(\eps x)}\,dx + O(\eps^\frac{n}{p}) \\
 & = \eps^p  f(0) \int_{\R^n} U(x)^p\,dx + o(\eps^p).
\end{split}
\end{equation*}

\subsection{Proof of \eqref{EstimGrad}}

We first write
\begin{equation*}
\begin{split}
 \int_{\R^n} f(x) |\nabla u_\eps|^{p(x)}\,dx
  =   \int_{\R^n} f(x) |\eta\nabla U_\eps + U_\eps\nabla \eta|^{p(x)}\,dx
  = \int_{\R^n} f(x) |\eta\nabla U_\eps|^{p(x)} \,dx + R_\eps,
\end{split}
\end{equation*}
where, using the inequality
$$ ||a+b|^q-|a|^q|\le C(|b|^q+|b||a|^{q-1}), $$
(the constant $C$ being uniform in $q$ for $q$ in a bounded interval of $[0,+\infty)$) we can estimate
$$ |R_\eps| \le C\Big[ \int_{B_{2\delta}\backslash B_{\delta}} |\nabla \eta|^{p(x)} U_\eps^{p(x)}\,dx
               + \int_{B_{2\delta}\backslash B_{\delta}}  |\nabla \eta| U_\eps(x) |\nabla U_\eps|^{p(x)-1}\,dx\Big]
              = C [   I_1(\eps) + I_2(\eps)]. $$
Since $U_\eps\le 1$ in $\R^n\backslash B_\delta$ for $\eps$ small, we can bound $I_1(\eps)$ as before by
\begin{equation*}
\begin{split}
 I_1(\eps)
 & \le C\int_{B_{2\delta}\backslash B_{\delta}} U_\eps^p \,dx
  \le C \eps^p \int_{\R^n\backslash B_{\delta/\eps}} U^p\,dx
  \le C \eps^p\eps^\frac{n-p^2}{p-1} = C \eps^\frac{n-p}{p-1},
\end{split}
\end{equation*}
if $p^2<n$.
Since $|\nabla U_\eps|\le 1$ in $\R^n\backslash B_\delta$ for $\eps$ small, we also have
\begin{equation*}
\begin{split}
 I_2(\eps)
 & \le C \int_{\R^n\backslash B_{\delta}} U_\eps(x)|\nabla U_\eps|^{p-1}\,dx\\
 & \le C \|U_\eps\|_{L^p(\R^n\backslash B_{\delta})}
       \|\nabla U_\eps\|^{p-1}_{L^p(\R^n\backslash B_{\delta})} \\
 & \le C \eps^\frac{n-p}{p(p-1)} \|\nabla U_\eps\|^{p-1}_{L^p(\R^n\backslash B_{\delta})},
\end{split}
\end{equation*}
with, since $ |U'(r)| \sim r^{-\frac{n-1}{p-1}} $ as $r\sim +\infty$,
\begin{equation*}
\begin{split}
 \int_{\R^n \backslash B_\delta} |\nabla U_\eps|^p\,dx
  & \le C \int_{\delta/\eps}^{+\infty} |U'(r)|^p r^{n-1}\,dr
  \le C \eps^\frac{n-p}{p-1}.
\end{split}
\end{equation*}
It follows that $I_2(\eps)=O(\eps^\frac{n-p}{p-1})$ and then $R_\eps = O(\eps^\frac{n-p}{p-1})$.
Independently, since
\begin{equation*}
\begin{split}
 |\nabla U_\eps (x)|
= & \frac{n-p}{p-1}\eps^{-n/p} \left( \frac{|x|}{\eps} \right)^\frac{1}{p-1}
  \left(1+\left(\frac{|x|}{\eps}\right)^\frac{p}{p-1}\right)^{-n/p} \end{split}
\end{equation*}
we have
\begin{equation}\label{condx}
 |\nabla U_\eps (x)| < 1 \quad \text{for }  |x|> C_p \eps^\frac{n-p}{p(n-1)},
               \quad C_p = \left(\frac{n-p}{p-1}\right)^\frac{p-1}{n-1}.
\end{equation}
Taking some constant $C>C_p$, we thus write
\begin{equation*}
\begin{split}
\int_{\R^n} f(x) |\nabla u_\eps|^{p(x)}\,dx = & \int_{B_{C\eps^\frac{n-p}{p(n-1)}}} f(x) |\nabla U_\eps|^{p(x)} \,dx\\
 &  + \int_{\R^n\backslash B_{C\eps^\frac{n-p}{p(n-1)}}} f(x) |\nabla U_\eps|^{p(x)} \,dx
   + O(\eps^{\frac{n-p}{p-1}}).
\end{split}
\end{equation*}
Since $|\nabla U_\eps(x)|<1$ in $ \R^n\backslash B_{C\eps^\frac{n-p}{p(n-1)}} $, we can  bound the second integral on the right hand side by
\begin{equation*}
\begin{split}
 C\int_{\R^n\backslash B_{C\eps^\frac{n-p}{p(n-1)}}} |\nabla U_\eps|^p \,dx
  & \le C  \int_{\eps^{-\frac{n(p-1)}{p(n-1)}}}^{+\infty} r^\frac{p}{p-1}
           \left( 1+r^\frac{p}{p-1} \right)^{-n}\,r^{n-1}\,dr
  \le C \eps^\frac{n(n-p)}{p(n-1)}
  = o(\eps^{\frac{n-p}{p-1}}).
\end{split}
\end{equation*}
Hence
\begin{equation*}
\begin{split}
 \int_{\R^n} f(x) |\nabla u_\eps|^{p(x)}\,dx
  =& \int_{B_{C\eps^\frac{n-p}{p(n-1)}}} f(x) |\nabla U_\eps|^{p(x)} \,dx  + O(\eps^\frac{n-p}{p-1})  \\
  =&  B_0(\eps)
    + B_1(\eps) \eps^2\ln\,\eps + o(\eps^2\ln\eps) + O(\eps^\frac{n-p}{p-1})\\
\end{split}
\end{equation*}
where $B_0(\eps)$ and $B_1(\eps)$ are the same as $B_0,B_1$ but integrating over $B_{\eps^{-\frac{n(p-1)}{p(n-1)}}}$ instead of $\R^n$. Again, as in the computation of  \eqref{EstimLq}, the term involving $\nabla f(0)$ vanishes for symmetry reasons.

Since $|U'(r)|^p\sim r^\frac{p(1-n)}{p-1}$ as $r\sim +\infty$, we have
$$ |B_0-B_0(\eps)|\le C\int_{\R^n\backslash B_{C\eps^{-\frac{n(p-1)}{p(n-1)}}}} |\nabla U|^p\,dx
  \le C \int_{\eps^{-\frac{n(p-1)}{p(n-1)}}}^{+\infty} r^{\frac{p-n}{p-1}-1}\,dr \le C \eps^\frac{n(n-p)}{p(n-1)}   = o(\eps^\frac{n-p}{p-1}), $$
$$ |B_1-B_1(\eps)|\le C\int_{\R^n\backslash B_{C\eps^{-\frac{n(p-1)}{p(n-1)}}}} |x|^2|\nabla U|^p\,dx
   \le C \eps^\frac{n(n-3p+2)}{p(n-1)} \quad \text{if } p<\frac{n+2}{3}.$$
Hence if $p<\frac{n+2}{3} $ we have
\begin{equation*}
 \int_{\R^n} f(x)|\nabla u_\eps|^{p(x)} \,dx - B_0  - B_1 \eps^2\ln\,\eps  = o(\eps^2\ln\eps).
\end{equation*}

\section*{Acknowledgements}
This work was partially supported by Universidad de Buenos Aires under grant X078 and by CONICET (Argentina) PIP 5478/1438.
A. Silva is a fellow of CONICET.

\bibliographystyle{plain}
\bibliography{biblio}

\end{document}